\documentclass[a4paper]{article}
\usepackage{amsmath}
\usepackage{amsfonts}
\usepackage{amssymb}         
\usepackage{graphicx}
\begin{document}
\DeclareGraphicsRule{.emf}{bmp}{}{}
\def\lag{\mathop{\rm lag}}
\def\det{\mathop{\rm det}}
\def\cov{\mathop{\rm cov}}
\def\Var{\mathop{\rm Var}}
\def\AIC{\mathop{\rm AIC}}
\def\Pr{\mathop{\rm Pr}}
\def\Re{\mathop{\rm Re}}
\def\Im{\mathop{\rm Im}}
\def\D{\mathop{\rm Var}}

\def\supp{\mathop{\rm supp}}
\newcommand{\E}{\mathbf E}
\def\lc{\mathop{\rm lc}}
\newcommand{\MISE}{\mathop{\rm MISE}}
\newcommand{\MSE}{\mathop{\rm MSE}}
\newcommand{\sqq}{\hbox{\vrule\vbox{\hrule\phantom{o}\hrule}\vrule}}
\newcommand{\mv}{\mathop{\rm v.p.}}
\begin{center}
{\large \textbf{ Rates of convergence for constrained deconvolution
 problem}}
\end{center}
\begin{center}
{\bf Denis Belomestny}
\end{center} 
\begin{center}
Institute for Applied Mathematics, University of Bonn, Wegelerstrasse 10, 53012, Bonn 
\end{center}
\centerline{\bf Abstract}
\begin{center}
{\small Let $X$ and $Y$ be two independent identically distributed random variables with density $p(x)$ and $Z=\alpha X+\beta Y$
for some constants $\alpha>0$ and $\beta>0$.
We consider the problem of estimating $p(x)$ by means of the samples from the distribution of $Z$. 
Non-parametric estimator based on the sync kernel is constructed and asymptotic behaviour of the
corresponding mean integrated square error is investigated.  
}
\end{center}
\section{Introduction}

Let $Z_1,...,Z_n$ be i.i.d observations, where $Z_i=X_i+Y_i$ and $X_i$ and $Y_i$ are independent random variables. Assume that 
the distribution of $Y$'s is known.    
The {\it ordinary deconvolution} problem is the problem of estimating  the destribution of a random variables $X_i$ from the observations $Z_i$.
In some situations (for instance, in signal processing) we  don't know the distribution of $Y_i$ but rather can assume  that $X_i$ and $Y_i$ belong to the same  class of
distributions ({\it constrained deconvolution}).
For example, $X_i$ and $Y_i$ may be of the same {\it multiplicative type} (see Belomestny (2002,2003)), that is
$$
Z_i=\alpha X_i+\beta Y_i,\quad \alpha>0,\,\beta>0,
$$
where $X_i$ and $Y_i$ are now independent identically distributed random variables.
Let distribution of $Z_i$ be absolutely continuous with  density $p_Z$ then the distribution of $X_i$(and $Y_i$) is also  continuous with some density, say $p(x)$.
Our aim  is to construct a non-parametric estimator for $p(x)$ based on the sample $Z_1,...,Z_n$  and to study its
asymptotic behaviour.
\\
{\bf Example}
Multiple Access FH SS Radio Networks systems have been considered for a variety of applications such as military ground-based communications
and cellular radio ( see Simon {\it et al} (1994) and Steele (1994)). In these wireless networks with the random multiple access protocol, because of channell reuse, each terminal interferes with signals
transmitted by other terminals. This interference is usually referred to as multiple access or self-interference.
In the system model , a receiver is located at the center of a plane where there are N  transmitters (terminals). The distance between the receiver and 
interfering terminals is denoted as $r_i$.The signal amplitude loss function over distance $r$ is given by 
$$
a(r)=\frac{K}{r^m},
$$   
where the constant $K$ depends on the transmitted power, and the attenuation factor $m$ characterizes the environment. The received passband
signal is
$$
Y(t)=\sum_{i=1}^N a(r_i)X_i(t),
$$
where $X_i(t)$ is the signal from the $i$th interferer.
Because all terminals use the same modulation scheme and power, it is reasonable to assume that $\{X_i(t)\}_{i=1}^N$ for every $t$ are independent and identically distributed (i.i.d.).
The problem in this case can be formulated as one of reconstructing the distribution of $X_i$ from the sampled distribution of $Y$. 
 \\
\section{Main results}
\begin{figure}[t]
\begin{center}
\includegraphics[width=6cm]{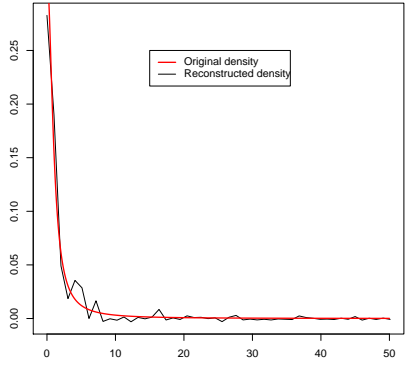}
\end{center}
\caption[]{\label{fig:compl}\small The results of the reconstruction  for
Cauchy density  ($n=1000$, $N=10$ , $\gamma=1/2$).}
\end{figure}

Let us denote by $f(t)$ the characteristic function of $X$ and suppose that $0<\alpha<\beta$. The characteristic function $g(t)$ of the random variable $Z/\beta$ can be expressed as
$$
g(t)=f(t)f(\gamma t),\quad 0<\gamma=\alpha/\beta<1.
$$
If  the infinite product  $\prod \frac{g(\gamma^{2k} t)}{g(\gamma^{2k+1}t)}$ converges  then
$$
f(t)=\prod_{k=0}^{\infty} \frac{g(\gamma^{2k} t)}{g(\gamma^{2k+1}t)}
$$
and  a natural estimator for $f(t)$ can be given as
\begin{equation}
\hat f_n(t)=\prod_{k=0}^{\infty} \frac{g_n(\gamma^{2k} t)}{g_n(\gamma^{2k+1}t)},
\label{eq:cf}
\end{equation}
provided that $g_n(t)\neq 0$ on $(0,t]$,where $g_n(t)$ is the empirical characteristic function corresponding to $g(t)$:
$$
g_n(t)=\frac{1}{n}\sum_{k=1}^n e^{it Z_k/\beta }.
$$
First of all we establish some asymptotic properties  of this estimator.
\\
{\bf Theorem 1 }\it
If  $\E|Z|^r<\infty$ for some $r>0$ and $g(u)\neq 0$ on $(0,t]$  then   
\begin{equation}
\sqrt{n}( \hat f_n(t) -f(t))\stackrel{D}{\longrightarrow}f(t)\sum_{k=0}^\infty \frac{(-1)^k}{g(\gamma^kt)}Y_F(\gamma^k t),
\label{mconv}
\end{equation}
where $Y_F=U(t)+iV(t)$ is  a complex valued
 Gaussian process with the $EY_F(t)=0$ and
with the  cross-covariance matrix  
$$
{\bf C}( t, s)=
\begin{pmatrix}
EU(t)U( s)&EU(t)V( s)\\
EV( t)U( s)&EV( t)V( s)\\
\end{pmatrix}
=
$$
$$
\begin{pmatrix}
{1\over2}[u( t-s)+u( t+s)]
-u(t)u(s)&{1\over2}[-v(t- s)+v( t+ s)]
-u(t)v( s)\\
{1\over2}[v( t- s)+v( t+s)]
-v(t)u(s)&{1\over2}[u( t- s)-u( t+ s)]
-v( t)v( s)\\
\end{pmatrix},
$$
where $u(t)$ and $v(t)$ denote real and imaginary part of $g(t)$.
\rm
\\
Turn now to the original problem of esimating  $p(x)$. 
Taking into account (\ref{eq:cf}),we can define the estimator for $p(x)$ as follows
\begin{equation}
\hat p_{nN}(x)=\frac{1}{2\pi}\int_{-1/h_n}^{1/h_n} e^{-itx} \left[\prod_{k=0}^{N} \frac{g_n(\gamma^{2k} t)}{g_n(\gamma^{2k+1}t)}\right]\, dt,
\label{eq:estor}
\end{equation}
with
\begin{equation}
1/h_n=\min \{\max \{\theta >0:|g_n(\theta )|>\varepsilon_n \},c_n\},
\label{eq:seq}
\end{equation}
where $\varepsilon_n>0$ and $c_n>0$ are two sequences of real numbers tending
correspondingly to zero and infinity as $n\to \infty $. The exact form of $\varepsilon_n$ and $c_n$ will be defined later. 
Let us  further put 
$$
d_n=\inf_{|s|<c_n}|g(s)|.
$$
In the following theorem we give some results concerning the behaviour of mean integrated square error corresponding to the estimator (\ref{eq:estor})
$$
\MISE(\hat p_{nN})=\int_{-\infty}^\infty |p(x)-\hat p_{nN}(x)|^2\,dx.
$$
\\
{\bf Theorem 2}
\it
Let the following conditions be satisfied
\begin{enumerate}
\item 
$$\int_{-\infty}^{\infty} |x|^rp_Z(x)\,dx<\infty,\quad  r>0$$
\item $p(x)\in L_2(\mathbb R)$
\item $\phi(t)\equiv |g(t)|$ does not increase on $(0,\infty)$
\end{enumerate}
There exist $A\equiv A(r)>0$ and  $D\equiv D(r)>0$  such that 
for
$$ 
\varepsilon_n=An^{-1/2}\log^{1/2}n,\quad c_n\leq\varphi^{-1}(2\varepsilon_n), 
$$ 
where $\phi^{-1}(\cdot)$ is inverse function to $\phi(\cdot)$
$$
\MISE(\hat p_{nN})\leq D\left[\frac{2^Nc_n}{n\varphi^2(c_n)}+\frac{ \gamma^{2rN}c_n^{r+1}}{(1-\gamma^{r/2})^2} \right]+\frac{1}{2\pi}\int_{|t|>c_n}|f(t)|^2\,dt.
$$
\rm
{\bf Corollary}
If $p_Z(x)$ is a  density of   stable distribution and
$$|\E e^{itZ}|=e^{-b |t|^{a}} ,\quad b>0,\quad 0<a\leq 2,$$
then $p(x)$ is also  a stable density and for the estimator $\hat p_{nN}$ with 
$$
N=  \nu \log n,\quad c_n=(\zeta \log n)^{1/a},\quad \nu>0, \quad \zeta>0
$$
\begin{multline}
\MISE(\hat p_{nN})\leq C(\log n)^{(1+a)/a}\times\\ 
\times n^{-1/(2+\gamma^{1/a}-\ln 2/((a-\delta)\ln\gamma))},\quad a>\delta>0.
\label{col}
\end{multline}
\\
Results of the simulation for Cauchy density are presented in Fig.\ref{fig:compl}. 
\\
{ \bf Remark 1 } As has been shown by Stefanski and Carroll (1990) the best possible rates of convergence for the MISE
in general deconvolution problem are usually (including normal and Cauchy distribution)  $\log^{-p}(n)$ for some $p>0$.
We see that in the case of the constrained deconvolution  the situation is  better.
\\
{\bf Remark 2 }
The above method  can  produce estimates which
are not probability density functions i. e. may take negative values 
or/and do not integrate to one.
It  happens due to the finiteness of  $n$. 
For this reason,  some methods of modification of 
density estimators (all estimators, not only kernel estimators) has been constructed (see Glad {\it et al} (1999,2003)) 
in such a way that the resulting estimator always 
produces estimates which are almost surely probability density
functions, and, in addition, the resulting estimator is better
or at least almost as good as the initial one.

\section{Auxiliary results}
{\bf Lemma 1} \it
For any $r\in(0,2]$ and any  real characteristic function $f(t)$ of 
a distribution with finite absolute moment $\beta_r$ of the order $r$
the following inequality holds
$$ f(t)\geq 1-c(r)\beta_r|t|^r,$$
where
\begin{equation}
c(r)=2/(2r)^{r/2}.
\label{cr}
\end{equation}
\rm
{\bf Proof.}
First of all we prove the inequality 
\begin{equation}
\cos(x)>1-c(r)x^r,\quad x>0.
\label{eq:in}
\end{equation}
From the equality
$$
\left(\frac{1-\cos(x)}{x^r}\right)'=\frac{x^r\sin x- r x^{r-1}(1-\cos(x))}{x^{2r}}=0
$$
we have
\begin{equation}
r(1-\cos(x))=x\sin(x),\quad \frac{1-\cos(x)}{x^r}\leq \frac{x^{2-r}}{r},\quad x>0 
\label{eq:in1}
\end{equation}
On the other hand
\begin{equation}
\frac{1-\cos(x)}{x^r}\leq \frac{2}{x^r},\quad x>0.
\label{eq:in2}
\end{equation}
Combining (\ref{eq:in1}) and (\ref{eq:in2}), we get (\ref{eq:in})
$$
\frac{1-\cos(x)}{x^r}\leq \frac{2}{(2r)^{r/2}},\quad x>0.
$$
Finally,
$$
f(t)=\int_{-\infty}^\infty \cos(tx)\,dF(x)\geq \int_{-\infty}^\infty (1-c(r)|tx|^r)\,dF(x)=1-c(r)\beta_r|t|^r.
$$
{\bf Lemma 2 }\it
Let $\xi$ be a random variable with characteristic function $f$ and finite absolute moment $\beta_r$ of order $r\in (0,2]$, then for
$|t|<1/2\sqrt[r]{10c(r)\beta_r}$ the following inequality holds
\begin{equation}
|\ln f(t)|\leq 2\pi\sqrt{2^{r}|t|^rc(r)\beta_r},
\label{eq20}
\end{equation}
where $c(r)$ is given by (\ref{cr}).
\rm
\\
{\bf Proof.}
Since $|f(t)|^2$ and $\Re^2 f(t)$ are two real c.f.  we have
$$
0\leq 1- |f(t)|^2\leq 2^{r+1}|t|^rc(r)\beta_r,\quad 0\leq 1-{\Re}^2 f(t)\leq 2^{r+1}|t|^rc(r)\beta_r
$$
and 
$$
\sin^2\psi(t)=\frac{\Im^2 f(t)}{|f(t)|^2}< \frac{1-|f(t)|^2+1-\Re^2 f(t)}{|f(t)|^2}\leq \frac{2^{r+2}|t|^rc\beta_r
}{1-2^{r+1}|t|^rc\beta_r},
$$
where $\psi(t)\equiv\arg f(t)$ satisfying $\psi(0)=0$.
The elementary inequality $|x|\leq \frac{\pi}{3} |\sin x|$ that holds for $|x|<\pi/6$ entails 
\begin{equation}
|\psi(t)|\leq \frac{\pi}{3}\sqrt{\frac{2^{r+2}|t|^rc\beta_r
}{1-2^{r+1}|t|^rc\beta_r}}
\label{eq21}
\end{equation}
if $|t|<1/2\sqrt[r]{10c\beta_r}$.
Further, using the inequality $|\log(1-x)|\leq 2|x|\,, 0<x\leq 1/2$
\begin{equation}
-\ln |f(t)|\leq -\frac{1}{2}\ln\left(1-2^{r+1}|t|^rc(r)\beta_r\right)\leq 2^{r+1}|t|^rc(r)\beta_r.
\label{eq22}
\end{equation}
The combination of (\ref{eq21}) and (\ref{eq22}) gives  (\ref{eq20}).
$$\eqno\square$$
\textbf{Lemma 3 } \it Let us have $n$ i.i.d random variables  $X_1,...,X_n$ with common characteristic function $f(t)$ and the underlying
distribution possesses finite absolute moment $\beta_r$ of order $r>0$ . Then,
for $a>0$, $0<b\leq 2$
\[
\Pr \left( \sup_{|\theta |<a}|f_n(\theta )-f(\theta )|>b\right)\leq 2(1+a\Theta(n,b,r)) e^{-nb^2/144}+\frac{\nu_r}{n},
\]
where  $f_n(\theta)$ is the corresponding empirical characteristic function
$$
f_n(t)=\frac{1}{n}\sum_{k=1}^n e^{itX_k},
$$
$\nu_r$ is a constant not depending  on $\beta_r$, $a$,$b$ and 
$$
\Theta(n,b,r)=\left\{
\begin{aligned}
\frac{\beta^{1/r}_rn^{(2-r)/r}}{b^{1/r}} &,& 0<r\leq1\\
\frac{\beta^{2/r}_{r/2}n^{1/r}}{b^{2/r}} &,& 1< r\leq 2.\\
\end{aligned}
\right.
$$
\rm
\textbf{Proof.}
Let us prove the inequality for $0<r\leq 1$.
 Define 
\[
\gamma =\left({\frac{{b}}{15{n^{2-r}\beta _r}}}\right)^{1/r}.
\]
We find numbers $t_1<t_2<...<t_k$ with the property that $t_1=-a$, $t_k=a$, $%
|t_i-t_{i+1}|\leq \gamma $. Clearly, we can assure this with $k\geq
1+2a/\gamma $. We begin with 
\[
\Pr \left( \sup_{|t|<a}|f_n(t)-f(t)|>b\right) \leq \Pr \left(
\sup_{|t-s|<\gamma }|f(t)-f(s)|>b/3\right) 
\]
\[
+\Pr \left( \sup_{|t-s|<\gamma }|f_n(t)-f_n(s)|>b/3\right) +\sum_{i=1}^k\Pr
(|f_n(t_i)-f(t_i)|>b/3).
\]
Denote the three summands on the right hand side by $T_1$, $T_2$ and $T_3$
and estimate them.  Due to the inequality $|1-e^{ix}|\leq c(r)|x|^r$ that holds for $0<r\leq 1$ with $c(r)= \sqrt{4/(2r)^r+1/r^{2r}}\leq 5$ we have
\[
|f(t)-f(s)|\leq E|1-e^{i(t-s)X}|\leq c(r)E|(t-s)X|^r\leq c(r)\gamma^r \beta _r\le b/3,
\]
when $|t-s|\leq \gamma $. Therefore, $T_1=0$. Next, we let $Y$ be the random
variable that puts mass $1/n$ at each of the $X_i$'s, then 
\[
|f_n(t)-f_n(s)|\leq E|1-e^{i(t-s)Y}|\leq E|(t-s)Y|=|t-s|\left| {\frac 1n}%
\sum_{i=1}^nX_i\right| .
\]
Using the inequality  
\begin{equation}
{\bf E}\left| {\frac 1n}\sum_{i=1}^nX_i\right|^r\leq n^{1-r} {\bf E}|X_1|^r, \quad 0<r\leq 1, 
\label{eq31}
\end{equation}
we  get
\[
T_2\leq \Pr \left( \gamma \left| {\frac 1n}\sum_{i=1}^nX_i\right| \ge
b/3\right)  \leq \frac{5b^{1-r}3^{r+1}}{n}\leq \frac{\nu_r}{n} 
\]
by the Chebyshev inequality. Finally, for fixed $t_i$, 
\[
\Pr (|f_n(t_i)-f(t_i)|>b/3)\leq \Pr (|u_n(t_i)-u(t_i)|>b/6)
\]
\[
+\Pr (|v_n(t_i)-v(t_i)|>b/6)\leq 2e^{-nb^2/144},
\]
by Berneistein's inequality for bounded random variables (see , for example, Bosq(1998)),
where as usually $u(t)$ and $v(t)$ are the real and imaginary parts of $f(t)$%
, and $u_n(t)$ and $v_n(t)$ are those of $f_n(t)$. 
The proof in the case $1< r\leq 2$ can be conducted in a similar way   using  the inequality
$$
{\bf E}\left| {\frac 1n}\sum_{i=1}^nX_i\right|^r\leq {\bf E}|X_1|^r, \quad r> 1
$$
instead of (\ref{eq31}).
 $$\eqno \square$$
\section{Proofs of the main results}
{\bf Proof of Theorem 1}
Let us consider the process 
$$Y_n(t)=\sqrt n(g_n(t)-g(t))$$
with $EY_n(t)=0$, the cross-covariance matrix ${\bf C}( t, s)$
and
$$EY_n(t)\overline{Y_n(s)}=f(t-s)-
f(t)f(-s).$$
The finite-dimensional distributions of $Y_n(t)$
converge by the multidimensional central limit theorem to
those of $Y_F(t)$ as $n\to\infty$ (Ushakov (1999), Chapter 3).
Further, it is known (see, for example,  Billingsley (1968),
Theorem 4.2) that if 
\begin{equation*}
\zeta_{un}\buildrel {D}\over\to \zeta_u\ \ {\rm as}\ \ n\to\infty
\end{equation*}
for each $u$, 
$$\zeta_u\buildrel {D}\over\to \zeta\ \ {\rm as}\ \ u\to\infty,$$
and 
\begin{equation}
\lim_{n\to\infty}\lim_{u\to\infty}
\sup\Pr(|\zeta_{un}-\eta_n|\ge\varepsilon)=0,
\label{conv1}
\end{equation}
then
$$\eta_n\buildrel {D}\over\to \zeta\ \ {\rm as}\ \ n\to\infty.$$
Put
$$\zeta_N=\sum_{k=0}^N \frac{(-1)^k}{g(\gamma^kt)}Y_F(\gamma^k t),\quad \zeta= \sum_{k=0}^\infty \frac{(-1)^k}{g(\gamma^kt)}Y_F(\gamma^k t)$$
$$\zeta_{Nn}= \sum_{k=0}^N (-1)^kY^L_n(\gamma^kt),\quad \eta_n= \sum_{k=0}^\infty (-1)^k Y^L_n(\gamma^k t),$$
where $Y^L_n=\sqrt{n}(\ln g_n(t)-\ln g( t))$.
We have
\begin{multline}
\Pr(|\zeta_{Nn}-\eta_n|\ge\varepsilon)=
\Pr\left(\left|\sum_{k=N+1}^\infty (-1)^kY^L_n(\gamma^k t)\right|\ge\varepsilon \left|\max_{k>N}\left| \frac{g(\gamma^kt)-g_n(\gamma^k t)}{g(\gamma^k t)}\right|\leq \frac{1}{2}\right.\right)+\\
\Pr\left(\max_{k>N}\left| \frac{g(\gamma^kt)-g_n(\gamma^k t)}{g(\gamma^k t)}\right|> \frac{1}{2}\right)=P_1+P_2
\label{eq:p}
\end{multline}
Using the Markov and Cauchy-Schwarz inequality, we have
$$P_1\le{1\over{\varepsilon^2}}\sum_{i=N+1}^\infty\sum_{j=N+1}^\infty (-1)^{i+j}{\bf E}\left\{Y^L_n(\gamma^it)\overline{Y^L_n(\gamma^jt)}\left|\max_{k>N}\left| \frac{g(\gamma^kt)-g_n(\gamma^k t)}{g(\gamma^k t)}\right|\leq \frac{1}{2}\right.\right\}\
$$
\begin{equation}
\le{1\over{\varepsilon^2}}\left(\sum_{j=N+1}^\infty \left[{\bf E}\left\{ \left|Y^L_n(\gamma^jt)\right|^2\left|\max_{k>N}\left| \frac{g(\gamma^kt)-g_n(\gamma^k t)}{g(\gamma^k t)}\right|\leq \frac{1}{2}\right.\right\}\right]^{1/2}\right)^2.
\label{eq:p1}
\end{equation}
Since $g(u)\neq 0$ for $|u|<t$ there exists $b>0$ such that $g(u)>2b$ on $(-t,t)$ and due to Lemma 3
\begin{multline}
 P_2\leq \Pr\left(\max_{k>N}|g_n(\gamma^kt)-g_n(\gamma^kt)|>b\right)\leq\\
 4\left(1+\gamma^{2N+1}|t|\Theta(n,b,r)\right) e^{\frac{-nb^2}{144}}+\frac{\nu_r}{n}.
\label{eq:p2}
\end{multline}
Further, using the elementary inequality $|\ln(1-z)|\leq 2|z|$ that holds for $|z|\leq1/2$ we get 
\begin{multline*}
{\bf E}\left\{\left |\ln \left(1-\frac{g(\gamma^kt)-g_n(\gamma^k t)}{g(\gamma^k t)} \right)\right|^2\left|\left| \frac{g(\gamma^kt)-g_n(\gamma^k t)}{g(\gamma^k t)}\right|\leq \frac{1}{2}\right.\right\}\leq\\
\leq \frac{4}{|g(\gamma^k t)|^2}{\bf E}|g(\gamma^kt)-g_n(\gamma^k t)|^2
=\frac{4}{n|g(\gamma^k t)|^2}(1-|g(\gamma^k t)|^2).
\end{multline*}
Combining (\ref{eq:p}),(\ref{eq:p1}), (\ref{eq:p2}) and using inequality $|g(t)|^2>1-c(r)\beta_r|t|^r$(see Lemma 1) we get (\ref{conv1}).
Analogously
$$
\Pr(|\zeta_N-\zeta|\ge\varepsilon)\le{1\over{\varepsilon^2}}
\sum_{i=N+1}^\infty\sum_{j=N+1}^\infty (-1)^{i+j}[g((\gamma^{j}-\gamma^{i})t)-g(\gamma^{i}t)g(-t\gamma^{j})]\leq
$$
$$
\sum_{i=N+1}^\infty\sum_{j=N+1}^\infty C[(\gamma^{rj}+\gamma^{ri})|t|+\gamma^{r(i+j)}|t|^2],\quad N>N_0
$$
that implies $\zeta_N\buildrel {D}\over\to \zeta$, $N\to\infty$.
Thus,
$$
\sqrt{n}(\ln \hat f_n(t)-\ln f(t))\stackrel{D}{\longrightarrow} \sum_{k=0}^\infty \frac{(-1)^k}{g(\gamma^k t)}Y_F(\gamma^k t)
$$
and  
$$
\sqrt{n}( \hat f_n(t) -f(t))\stackrel{D}{\longrightarrow}f(t)\sum_{k=0}^\infty \frac{(-1)^k}{g(\gamma_kt)}Y_F(\gamma^k t)
$$
as $n\to\infty$.
$$\eqno\square$$
{\bf Proof of Theorem 2}
Plancherel-Parseval formula entails
\begin{multline*}
\int_{-\infty}^\infty |p(x)-\hat p_{nN}(x)|^2\,dx= \frac{1}{2\pi}\int_{-1/h_n}^{1/h_n}  \left|\prod_{k=0}^{N} \frac{g_n(\gamma^{2k} t)}{g_n(\gamma^{2k+1}t)}- \prod_{k=0}^{\infty}\frac{g(\gamma^{2k} t)}{g(\gamma^{2k+1}t)}\right|^2\, dt
\\
+\frac{1}{2\pi}\int_{|t|>1/h_n}|f(t)|^2\,dt.
\end{multline*}
Let us now estimate the first summand on the right-side hand.
Using lemma 2 we have for some constants $C\equiv C(r)$ and $D\equiv D(r)$ not depending on $N$ and $\gamma$
\begin{multline*}
\left|\ln\left(\prod_{k=N+1}^\infty\frac{g(\gamma^{2k} t)}{g(\gamma^{2k+1}t)}\right)\right|\leq\sum_{k=2(N+1)}^\infty|\ln g(\gamma^{k}t)|
\leq C\frac{ \gamma^{r(N+1)}|t|^{r/2}}{1-\gamma^{r/2}}
,\quad |t|<\frac{D}{\gamma^{2(N+1)}}
\end{multline*}
 The elementary inequality $|z|\leq2|\ln(1-z)|$ that holds for $|z|\leq 1/3$ implies 
$$
\left|1-\prod_{k=N+1}^\infty\frac{g(\gamma^{2k} t)}{g(\gamma^{2k+1}t)}\right|\leq 
C\frac{ \gamma^{r(N+1)}|t|^{r/2}}{1-\gamma^{r/2}} ,\quad |t|<\frac{D}{\gamma^{2(N+1)}}
$$
Further,
\begin{multline*}
\prod_{k=0}^N \frac{g_n(\gamma^{2k} t)}{g_n(\gamma^{2k+1}t)}- \prod_{k=0}^N\frac{g(\gamma^{2k} t)}{g(\gamma^{2k+1}t)}=
\sum_{i=0}^N H(\gamma^{2i} t) \prod_{k\neq i} \frac{g(\gamma^{2k} t)}{g(\gamma^{2k+1}t)}+\\
 \sum_{i_1\neq i_2} H(\gamma^{2i_1} t)H(\gamma^{2i_2} t)\prod_{k\neq i_1,i_2}\frac{g(\gamma^{2k} t)}{g(\gamma^{2k+1}t)}+\ldots+\prod_{i=0}^N H(\gamma^{2i} t),
\end{multline*}
where
$$
H(t)=\frac{g_n(t)}{g_n(\gamma t)}- \frac{g( t)}{g(\gamma t)}.
$$
Since $|g(t)|$ is non-increasing
$$
\left|\frac{g(t)}{g(\gamma t)}\right|\leq 1
$$
and for $|t|<1/h_n$
$$
|H(t)|\leq\frac{\Delta(t)}{\varepsilon_n},
$$
where
$$
\Delta(t)=|g_n(t)-g(t)|+|g_n(\gamma t)-g(\gamma t)|.
$$
Further, using the equality $\E|g_n(t)-g(t)|^2=\frac{1-|g(t)|^2}{n}$, one gets
\begin{multline*}
\E\int_{-1/h_n}^{1/h_n}  \left|\prod_{k=1}^{N} \frac{g_n(\gamma^{k-1} t)}{g_n(\gamma^{k}t)}- \prod_{k=1}^{N}\frac{g(\gamma^{k-1} t)}{g(\gamma^{k}t)}\right|^2\, dt\leq 
4c_n\Pr(\sup_{|t|<c_n}\Delta(t)>\varepsilon_n)+\frac{4^Nc_n}{n\varepsilon_n^2}.
\end{multline*}
The combination of the previous estimates and some simple calculations yield
\begin{multline*}
\E\int_{-1/h_n}^{1/h_n}  \left|\prod_{k=1}^{N} \frac{g_n(\gamma^{k-1} t)}{g_n(\gamma^{k}t)}- \prod_{k=1}^{\infty}\frac{g(\gamma^{k-1} t)}{g(\gamma^{k}t)}\right|^2\, dt\leq 
8c_n\Pr(\sup_{|t|<c_n}\Delta(t)>\varepsilon_n)+\frac{4^{N+1}c_n}{n\varepsilon_n^2}+\\
C\frac{ \gamma^{2rN}c_n^{r+1}}{(1-\gamma^{r/2})^2} , \quad C>0
\end{multline*} 
where according to lemma 3
$$
\Pr(\sup_{|t|<c_n}\Delta(t)>\varepsilon_n)\leq \Pr(\sup_{|t|<c_n}|g_n(t)-g(t)|> \varepsilon_n/2)+\Pr(\sup_{|t|<c_n}|g_n(\gamma t)-g(\gamma t)|> \varepsilon_n/2)\leq 
$$
$$
4\left(1+c_n\Theta(n,\varepsilon_n/2,r)\right) e^{\frac{-n\varepsilon_n^2}{4\cdot144}}+\frac{\nu_r}{n}
$$
The condition $p(x)\in L_2(\mathbb R)$ implies $\phi(\cdot)\in L_2(\mathbb R)$ that in its turn means
$$
\phi^{-1}(x)<\frac{1}{x^2} ,\quad x<\delta
$$
and therefore $c_n=O (n)$ as $n\to \infty$.
So. one has for some $\kappa=\kappa(r)>0$ and $C>0$
$$
\Pr(\sup_{|t|<c_n}\Delta(t)>\varepsilon_n)\leq 4\left(1+\frac{C}{A}n^{1+\kappa}\right) n^{-4A/144}=O(1/n),\quad n\to \infty
$$
for large enough $A$. 

Thus,there exist constant  $D>0$ such that
\begin{equation}
\int_{-\infty}^\infty |p(x)-\hat p_{nN}(x)|^2\,dx\leq D\left[\frac{4^Nc_n}{nd_n^2}+\frac{ \gamma^{2rN}c_n^{r+1}}{(1-\gamma^{r/2})^2} \right]+\frac{1}{2\pi}\int_{|t|>1/h_n}|f(t)|^2\,dt.
\label{eq:li}
\end{equation}
Let us now estimate the last term in (\ref{eq:li}).
We have
\begin{equation}
\mathbf{E}\left[\int\limits_{|t|\ge 1/h_n}|f(\theta )|^2d\theta\right] \leq \Pr \left[1/h_n<c_n\right]\int\limits_{\mathbb R}|
f(\theta )|^2d\theta + \int\limits_{|\theta |\ge c_n}|f(\theta )|^2d\theta.
\label{ineq1}
\end{equation}
Since 
$$
d_n-\varepsilon_n>\phi(c_n)-\varepsilon_n\geq \varepsilon_n,
$$
lemma 3 yields
\begin{multline*}
\Pr \left [1/h_n<c_n\right ]=\Pr (\inf_{|\theta |<c_n}|g_n(\theta )|\leq \varepsilon_n )\leq \\
\Pr (\sup_{|\theta |<c_n}|g_n(\theta )-g(\theta )|>d_n-\varepsilon_n )\leq  \\
2\left( 1+c_n\Theta(n,\varepsilon_n,r)\right)
e^{-n\varepsilon^2_n /144}+\frac{\nu_r}{n}=O(1/n),\quad n\to \infty
\end{multline*}
$$\eqno\square$$
{\bf Proof of Corollary 1.}
Let us put
$$
N=  \nu c^\alpha_n,\quad c_n=(\zeta \ln n)^{1/\alpha},\quad \zeta>0 
$$
Since $X$ posesses moments of $r$ order for $r<\alpha$  
$$
\int_{-\infty}^\infty |p(x)-\hat p_{nN}(x)|^2\,dx\leq D\left[(\zeta\ln n)^{1/\alpha} n^{2\zeta \beta(1+\gamma^{1/\alpha})+\nu \zeta\ln 2 -1}+(\zeta\ln n)^{(r+1)/\alpha}n^{r\nu \zeta\ln\gamma}  \right]
$$
$$
+(\ln n)^{2\alpha}n^{-2\zeta\beta}
$$
$$
$$
Taking $\nu=2\beta/r\ln(1/\gamma)$ and $\zeta=1/2\beta(1+(1+\gamma^{1/\alpha})-\ln 2/r\ln\gamma)$ we come to (\ref{col}).

\end{document}